\documentclass[12pt]{paper}

\usepackage{amscd}
\usepackage[titletoc,toc,title]{appendix}

\usepackage{amsthm}
\usepackage{thmtools}

\usepackage[headings]{fullpage}
\usepackage{mathrsfs}
\usepackage{amsmath,amsfonts,amssymb,amsthm}
\usepackage[all]{xy}
\usepackage{amscd}

\usepackage[dvips]{graphics,color}
\usepackage{hyperref}

\def\noi{\noindent}

\def\pf{\noi{\bf Proof.\ \,}}

\def\eop{{$\square$}}

\def\QQ{{\mathbb Q}}

\def\ZZ{{\mathbb Z}}

\def\mm{{\mathfrak m}}

\def\mp{{\mathfrak p}}

\def\Supp{{\mathrm{Supp}}}
\def\Ass{{\mathrm{Ass}}}
\def\Ker{{\mathrm{Ker}}}
\def\Im{{\mathrm{Im}}}

\def\Spec{{\mathrm{Spec}}}

\def\Supp{{\mathrm{Supp}}}

\def\Hom{{\mathrm{Hom}}}

\def\End{{\mathrm{End}}}

\def\ASpec{{\mathrm{ASpec}}}
\def\ASupp{{\mathrm{ASupp}}}

\def\sub{{\mathrm{sub}}}
\def\quot{{\mathrm{quot}}}
\def\ext{{\mathrm{ext}}}

\def\noeth{{\mathrm{noeth}}}

\def\al{ {\cal A}  }
\def\bl{ {\cal B}  }

\def\fl{ {\cal F}  }

\def\nl{ {\cal N}  }

\def\tl{ {\cal T}  }
\def\sl{ {\cal S}  }

\begin{document}

\newtheorem{thm}{Theorem}[section]
\newtheorem{prop}[thm]{Proposition}
\newtheorem{lem}[thm]{Lemma}
\newtheorem{coro}[thm]{Corollary}

\theoremstyle{definition}
\newtheorem{de}[thm]{Definition}
\newtheorem{nota}[thm]{Notation}

\newtheorem{rem}[thm]{Remark}
\newtheorem{ex}[thm]{Example}

\begin{center}

{\bf \ A classification of nullity classes in abelian categories}
\medskip

Yong Liu \footnote{Email address is \emph{yongliue@gmail.com}} and
Don Stanley\footnote{Email address is
\emph{donald.stanley@uregina.ca}}

Department of Mathematics and Statistics

University of Regina

\end{center}

\begin{abstract}
We give a classification of nullity classes (or torsion classes) in an abelian category by forming a spectrum of equivalence classes of premonoform objects. This is parallel to Kanda's classification of Serre subcategories.

%For an essentially small category $\al$ and its set of certain subcategories $\Phi$ with inclusion as partial order, we construct several topological spaces $K(\Phi)$, $K_p(\Phi)$ and $K_{gp}(\Phi)$ etc., called classifying space, so that any topological space that classifies these subcategories has $K_p(\Phi)$ as its Kolmogorov quotient up to homeomorphism. This assertion respects most of the well-known results such as Neeman's localizing subcategories in $D(R)$, Balmer's thick tensor ideals in a tensor triangulated category and Bensen's localizing tensor ideals in a homotopy category etc. Furthermore, we use our classifying space $K_{gp}(\Phi)$ to give a classification of the subcategories involved, which also recovers some of the famous results.

%Classifying space $K(\Phi)$ is defined from a given collection $\Phi$ of subcategories with certain closure operations, in which the primeness (or irreducibility) of subcategories plays a crucial role. We have checked in some well-know examples that thick subcategories of compact objects and localizing subcategories in an appropriate Noethrian stable homotopy category are primely (or irreducibly) generated respectively. For those with implicit generators, we produce homeomorphisms between our classifying space to the known spaces with Hochster dual topology, such as Balmer's spectrum of a tensor triangulated category and a quasi-compact, quasi-separated scheme.
\end{abstract}

\tableofcontents

\section{Introduction}

The classification of various types of subcategories is a fundamental problem in many different contexts. For example, Gabriel's thesis~\cite{Gabriel} showed that in the category of modules over a commutative noetherian ring $R$, the Serre subcategories are classified via the specialization closed subsets of the prime spectrum $\Spec~R$. Aisles, which correspond to $t$-structures~\cite{BBD}, have also been classified in various contexts such as~\cite{Stan} and~\cite{Alonso2}. There are many other famous results proving that the target subcategories are one to one correspondent to the closed (or dually open) subsets of certain spectrum, such as Neeman's classification of (co)localizing subcategories in~\cite{Neeman} and~\cite{Neeman2}, Balmer's classification of the radical thick tensor ideals in a tensor triangulated cateogry in~\cite{Balmer}, and so on. The classification of subcategories can also be used for reconstruction, especially for schemes. This is achieved by Rosenberg~\cite{Rosenberg} and Rouquier~\cite{Rouquier}, by using appropriate spectra.

More recently, R. Kanda in his paper~\cite{Kanda} classifies the Serre subcategories in an abelian category using the atom spectrum consisting of equivalence classes of monoform objects (see Definition~\ref{monoformobj}). We approach our problem of classifying nullity classes (or torsion classes, see Definition~\ref{nullityclass}) by constructing a new spectrum $\Spec\al$ introduced in Section~\ref{spectrum}, which consists of equivalence classes of premonoform objects, a generalization of monoform objects. An object is \textit{premonoform} if all its endomorpisms are trivial or injective (see Definition~\ref{premonoform} and Lemma~\ref{predef}). As pointed out to us by M. Reyes, these premonoform objects were studied by Tiwary and Pandeya~\cite{TiwaryPandeya} in the category of modules over a noncommutative ring, and interestingly its dual notion ``copremonoform'' was studied by Xue~\cite{Xue}. Using premonoform objects, a parallel result to Kanda's is obtained.

%Although this definition seems very natural we were unable to find anywhere this class of objects had been previously studied.

\begin{thm}
For an abelian category $\al$, there is an order preserving bijection
\[
\Supp:\{\text{Nullity classes of noetherian objects}\}\stackrel{\sim}{\rightleftarrows}~~~~~~~~~~~~~~~~~~~~~~~~~~~~~~~~~~~~~~~~~\]\[~~~~~~~~~~~~~~~~~~\{\text{Closed and extension closed subsets of}~~\Spec_{\noeth}\al\}:\Supp^{-1},
\]
where $\Spec_{\noeth}\al$ is the subcollection of points in $\Spec\al$ such that every equivalence class has a noetherian representative.
\end{thm}

This allows us to prove the main theorem of classification of torsion classes.
\begin{thm}
For a noetherian abelian category $\al$, there is an order preserving bijection
\[
\Supp:\{\text{Torsion classes in}~\al\}\stackrel{\sim}{\rightleftarrows}\{\text{Closed and extension closed subsets of}~\Spec\al\}:\Supp^{-1}.
\]
\end{thm}
%As other classifications of various subcategories establish a correspondence from either closed or open subsets of a topological space to the target subcategories, we mention that the nullity classes cannot be classified by a topological space (see Definition~\ref{classifies}) but has to correspond to the closed subsets with one extra condition extension-closed as Theorem 1 claims, see also Section~\ref{exceptional}.

The paper is organized as follows. In Section~\ref{monopre} and Section~\ref{serrenull}, we introduce and compare two pairs of concepts in an abelian category, monoform and premonoform objects, Serre subcategories and nullity classes. In Section~\ref{relatedcon}, we use examples to give a comparison of (pre)monoform with related concepts. Later in Section~\ref{spectrum}, a new spectrum of an abelian category is defined, together with a notion of support. In Section~\ref{maintheorem} and~\ref{variants}, we finally give a classification of nullity classes and torsion classes.

\section{Monoform and premonoform object}\label{monopre}

In this section we introduce the premonoform  and monoform objects in an abelian category. The former gives the points in our spectrum.

\begin{de}\label{monoformobj}
A nonzero object $M$ in an abelian category $\al$ is \textit{monoform} if for every nonzero subobject $N\leq M$ there are no nonzero common subobjects of both $M$ and $M/N$.
\end{de}

The following lemma gives an equivalent description of monoform object which appeared as the definition of monoform module right after Corollary 2.4 in~\cite{Gordon}.

\begin{lem}
Let $\al$ be an abelian category and $M\in\al$ a nonzero object. Then $M$ is monoform if and only if for every subobject $N\leq M$, every $f:N\rightarrow M$ is either zero or monic.
\end{lem}

\pf Suppose there is a subobject $N$ and map $f:N\rightarrow M$ with $\Ker(f)\neq0$. Then the induced map $N/\Ker(f)\hookrightarrow M$ shows that $M$ is not monoform. Conversely, suppose $M$ is not monoform with quotient $N/N'$ a subobject. Then we have a map $N\twoheadrightarrow N/N'\hookrightarrow M$ which is neither zero nor monic.~~~\eop

%It is straight forward to show that in the context of a category of modules over a ring, the monoform object coincides with the \textit{strongly uniform module} as it is called in~\cite{Storrer}, see also Section 8 in~\cite{lam}. In particular,

For a commutative ring $R$ and any prime ideal $\mp$ of $R$, the quotients $R/\mp$ are monoform. This follows for example from Lemma 1.5 in~\cite{Storrer}. For our classification, we introduce the concept of premonoform.

\begin{de}\label{premonoform}
An object $M\neq0$ in an abelian category $\al$ is \textit{premonoform} if it contains no nontrivial quotient as a subobject, meaning that there is no injection from $M/N$ to $M$, for any nonzero proper subobject $N\leq M$.
\end{de}

It follows directly from the definition that monoform objects are premonoform.

Other than Gordon and Robson~\cite{Gordon}, the other early references that we found to the concept of monoform object, in the category of modules over a (non)commutative ring, is Goldman~\cite{Goldman} and Storrer~\cite{Storrer}. We also refer to Kanda~\cite{Kanda}, for defining monoform object in an abelian category, and Lam~\cite{lam} for related concepts in the module theory, such as essential or rational extensions.

\begin{lem}\label{predef}
Let $\al$ be an abelian category and $M\in\al$ a nonzero object. Then the following are equivalent:

\noindent(1) $M$ is premonoform;

\noindent(2) $\Hom_\al(M/N,M)=0$ for any nonzero proper subobject $N\subseteq M$;

\noindent(3) any $f\in \End_\al M$ is either zero or injective;

\noindent(4) $\End_\al M$ is a domain (i.e. a ring with no zero divisors) and $M$ is a torsion free left module over $\End_\al M$.
\end{lem}

\pf The implication $(2) \Rightarrow (1)$ holds since otherwise there is a nontrivial quotient embedding into $M$. For $(1)\Rightarrow(2)$, suppose there is a nonzero map $f:M/N\rightarrow M$ for some nonzero proper subobject $N\leq M$. Then by the First Isomorphism Theorem there is some object $A$ with $N\leq A\leq M$ and $\Ker f=A/N$ such that $M/A\cong (M/N)/(A/N)\cong\Im f\hookrightarrow M$.

Let $M$ be premonoform and $0\neq f\in\End_\al M$. Then $N=\ker f$ is a nonzero subobject of $M$ such that $f$ induces an injection $M/N\cong \Im f\hookrightarrow M$. Hence $(1)\Rightarrow(3)$. Conversely, suppose there is an injection $M/N\hookrightarrow M$ for some subobject $N\leq M$. Then the composition $f:M\twoheadrightarrow M/N\hookrightarrow M$ as an element in $\End_\al M$ is either zero or an injection, that is, either $N=M$ or $N=0$. Therefore, $M$ is premonoform.

The equivalence $(3)\Leftrightarrow(4)$ was observed by Reyes (see his proof of Proposition 2.5 in~\cite{Reyes}). For $(3)\Rightarrow(4)$, suppose $f,g\in\End_\al M$ such that $gf=0$. Then $\Im f\subseteq\Ker g$. By assumption either $g=0$ or $\Ker g=0$, the latter of which implies $f=0$, as required. For $(4)\Rightarrow(3)$, take a nonzero $f\in \End_\al M$. Suppose $f(x)=0$ for $x\in M$. Since $M$ is torsion free, we have $x=0$. Thus $f$ is injective.~~~\eop

In Section~\ref{relatedcon} we will compare premonoform and monoform in more detail and in particular show (Proposition~\ref{precoincide}) that in the category of finitely generated modules over a commutative noetherian ring the two concepts coincide.

%Therefore, Section~\ref{monoformsect} provides lots of examples of premonoform objects as well.

\section{Serre subcategory, nullity class and torsion theory}\label{serrenull}

In this section, we discuss Serre subcategories, nullity classes and a way of constructing them. See~\cite{Stanw} for a discussion of nullity class in the category of modules over a ring. By absorbing objects which are isomorphic to one in the subcategory, we will thus assume all subcategories are \textit{replete} (that is, closed under isomorphisms) throughout the paper.

%We will assume our (sub)categories are \textit{replete} (that is, closed under isomorphisms) throughout the paper, and particularly any subcategory $\cl$ can be made replete by absorbing all the objects which are isomorphic to one in $\cl$.

%Recall that a subcategory in an abelian category is called a \textit{Serre} subcategory if it is a full abelian subcategory closed under taking subobjects, quotients and extensions. It is natural to use monoform objects to give a classification of Serre subcategories, see~\cite{Kanda}, whose definition relies on the construction of these categories.

%Now we extend the notion of monoform to study the nullity class in an abelian category.

\begin{de}\label{nullityclass}
Let $\al$ be an abelian category and $\bl$ a full subcategory of $\al$. Consider

%Then $\nl$ is called a \textit{nullity class} if

\noindent(1) for any exact sequence $0\rightarrow N\rightarrow M$, $M\in\bl$ implies $N\in\bl$;

\noindent(1') dually for any exact sequence $M\rightarrow N'\rightarrow0$, $M\in\bl$ implies $N'\in\bl$;

\noindent(2) for any exact sequence $0\rightarrow N\rightarrow M\rightarrow N'\rightarrow 0$, that $N, N'\in \bl$ implies $M\in\bl$,

\noindent then $\bl$ is called a \textit{nullity class} if it satisfies (1') and (2), and $\bl$ is called a \textit{Serre subcategory} if it is an abelian subcategory of $\al$ and satisfies (1), (1') and (2).

\end{de}

Notice that a nullity class $\nl$ is closed under \textit{retracts}, that is, whenever $A\oplus B\in\nl$, we have $A\in\nl$, thanks to the projection $p_A:A\oplus B\rightarrow A$.

A full subcategory $\tl$ of $\al$ is a \textit{torsion class} if it is both a nullity class and a coreflective subcategory, i.e. the inclusion functor $\tl\hookrightarrow\al$ admits a right adjoint. See Proposition 1.2 in~\cite{BeligReiten}. Furthermore, a Serre subcategory $\sl$ is \textit{localizing} if the quotient functor $\al\rightarrow\al/\sl$ admits a right adjoint. See Section A.2 in~\cite{Neeman3} for example. Next we look at closing collections of objects under the operations in the last definition, in order to get the Serre and nullity classes generated by that collection. We do this instead of taking the intersection of all Serre subcategories (or nullity classes) containing the collection, since this construction method will be used in later proofs.

\begin{nota}\label{next}
For every collection $\sl$ of objects in an abelian category $\al$, denote by
\[
\langle \sl\rangle_{\sub}=\{X\in\al~|~X~\text{is a subobject of some}~M\in\sl\}
\]
a full subcategory of $\al$ that is closed under subobjects. Dually, denote by $\langle \sl\rangle_{\quot}$ a full subcategory of $\al$ that is closed under quotients from $\sl$. Let $\sl^0$ be the class of the zero object. For any collections $\sl,\sl'$ of objects in $\al$, denote by
\[
\sl\ast\sl'=\{E\in\al~|~0\rightarrow X\rightarrow E\rightarrow X'\rightarrow 0,~\text{for some}~X\in\sl ,X'\in\sl'\}
\]
a full subcategory consisting of objects obtained by extensions from $\sl$ and $\sl'$. Let $\sl^1=\sl$. We define the $n$-\textit{extension}  $\sl^n=\sl\ast\sl^{n-1}$ recursively and denote
\[
\langle\sl\rangle_{\ext}=\bigcup_{n\geq 0}\sl^n.
\]
\end{nota}

The subcategories can thus be constructed via the above operations.

\begin{lem}\label{ext}
Let $\al$ be an abelian category and $\sl$ a collection of objects in $\al$. Then

%\noindent(1) the operation $\ast$ on collections of objects is associative;

%\noindent(2) $\langle M\ast N\rangle_{\quot}\subseteq\langle M\rangle_{\quot}\ast\langle N\rangle_{\quot}$;

\noindent(1) the category $\langle\langle\sl\rangle_{\quot}\rangle_{\ext}$ is a nullity class;

\noindent(2) the category $\langle\langle\langle\sl\rangle_\sub\rangle_\quot\rangle_\ext$ is a Serre subcategory.

\end{lem}

\pf Since $\langle\langle\sl\rangle_{\quot}\rangle_{\ext}$ is closed under taking quotient, (1) holds. See Proposition 2.4 in~\cite{Kanda} for a proof of (2).~~~\eop

%Both properties (1) and (2) hold thanks to the Snake Lemma. For a proof see (2) and (4) in Proposition 2.4 in~\cite{Kanda}. For (3), it is sufficient to show that this category $\langle\langle\sl\rangle_{\quot}\rangle_{\ext}$ is closed under taking quotients, and the result follows immediately from (2).~~~\eop

Nullity classes are torsion classes in many well-known abelian categories. First, recall that an abelian category $\al$ is \textit{noetherian} if it is essentially small and every object $A$ in $\al$ is \textit{noetherian} (that is, every ascending chain of the subobjects of $A$ becomes stationary after finitely many stages).

Another related set of conditions requires that the abelian category $\al$ is \textit{well-powered} (i.e. the collection $\mathrm{Sub}(A)$ of subobjects of any given object $A\in\al$ forms a set) and subcomplete in Dickson's sense~\cite{Dickson}, so that each $\mathrm{Sub}(A)$ forms a complete lattice by Proposition 1.1 in~\cite{Dickson}. The category $\al$ is \textit{subcomplete} if for every $A\in\al$, any subfamily $\{A_i\}_{i\in I}$ of $\mathrm{Sub}(A)$ has a sum $\Sigma_{i\in I}A_i$ and a product $\prod_{i\in I}A/A_i$. We thus have supremum and infimum defined by
\[
\bigcup_{i\in I}A_i=\Im(\Sigma_{i\in I}A_i\rightarrow A),~\bigcap_{i\in I}A_i=\Ker(\prod_{i\in I}A/A_i\rightarrow A)
\]
respectively. Note that any essentially small cocomplete abelian category is also well-powered and subcomplete.

%For example, if in the category every collection of subobjects of any given objects forms a complete lattice. See Section 4.2 in~\cite{Borceux}.

\begin{lem}\label{noethprop}
Every subset of a subobject lattice for a fixed object in a noetherian abelian category contains a maximal element.
\end{lem}

\pf Let $S$ be a subset of the suboject lattice, and $X_1\in S$. Suppose $S$ has no maximal element, then inductively there is $X_n\in S$ not contained in $\cup^{n-1}_{i=1}X_i$. We thus obtain a strictly increasing sequence $X_1\lneqq X_1\cup X_2\lneqq X_1\cup X_2\cup X_3\lneqq\cdots$ of infinitely many subobjects, a contradiction.~~~\eop

\begin{prop}\label{torsionnullity}
Let $\al$ be an abelian category, and $\tl$ a full subcategory in $\al$.

\noindent(1) If $\al$ is noetherian, then $\tl$ is a torsion class if and only if it is a nullity class;

\noindent(2) If $\al$ is well-powered and subcomplete, then $\tl$ is a torsion class if and only if it is a nullity class that is also closed under arbitrary coproducts.
\end{prop}

\pf The proof of (1) is similar to that of Theorem 2.1 in~\cite{Dickson}. Suppose $\al$ is noetherian and $\tl$ is a nullity class. Define as a full subcategory of $\al$ that
\[
\fl=\{A\in\al~|~\Hom_\al(X,A)=0~\text{for every}~X\in\tl\}.
\]
Then $\Hom_{\al}(X,A)=0$ for every $X\in\tl$ and every $A\in\fl$. Since $\al$ is noetherian, for every object $X\in\al$ there is a maximal subobject $X_\tl$ of $X$ among all subobjects of $X$ by Lemma~\ref{noethprop}. Then we claim that $X/X_\tl\in \fl$ by showing that any morphism $f:Y\rightarrow X/X_\tl$ is zero for every $Y\in\tl$, thus $(\tl,\fl)$ forms a torsion pair. Indeed, the image has the form $f(Y)=X'/X_\tl$ which lies in $\tl$, for some subobject $X'\subseteq X$, so that the short exact sequence
\[
0\rightarrow X_\tl\rightarrow X'\rightarrow X'/X_\tl\rightarrow0
\]
implies $X'\in\tl$ as well. Therefore, $X'=X_\tl$ and thus $f=0$, by the maximality of $X_\tl$.

For (2), we refer to Theorem 2.3 in~\cite{Dickson}.~~~\eop

Although the last result shows that nullity classes are often torsion classes, they are not always the same. For example, in the category of abelian groups, its subcategory of finitely generated abelian groups is a nullity class but not a torsion class.

\section{Comparing (pre)monoform and related concepts}\label{relatedcon}

In this section, we compare monoform and premonoform objects in different abelian categories with some other interesting concepts.

\begin{itemize}
\item \textit{General abelian category}
\end{itemize}

We have already seen that in a general abelian category, monoform implies premonoform by the definition. It is easy to see that any simple object (that is, it has no nonzero subobject) is monoform, also premonoform is indecomposable (see Proposition~\ref{relations}). Thus we obtain
\[
\text{simple}\Rightarrow\text{monoform}\Rightarrow\text{premonoform}\Rightarrow\text{indecomposable}.
\]
We will see in the next subsection and also Example~\ref{compareex} that all these implications are strict.

\begin{itemize}
\item \textit{A premonoform but not monoform object}
\end{itemize}

It is interesting to find in an abelian category a premonoform object that is not monoform. Let $T(V)$ be the tensor algebra of a vector space $V$ over a field $k$, with the basis $\{a,b\}$. Consider $T(V)$ as a left module over itself. It is known that $T(V)=\bigoplus_{i=0}^\infty V^{\otimes i}$ is a free associative algebra on the generators $a,b$. In this section, we show that $T(V)$ is premonoform but not monoform in the category of left $T(V)$-modules. Recall that in general if a vector space $V$ has a basis $\{e_1,...,e_m\}$, then the set of elements $e_s=e_{i_1}e_{i_2}\cdots e_{i_n}=e_{i_1}\otimes e_{i_2}\otimes\cdots \otimes e_{i_n}$ for $1\leq i_1,...,i_n\leq m$ with the unit 1 of $k$ forms a basis of $T(V)$, where $s=(i_1,...,i_n)$ and $n\in\ZZ^+$. The multiplication table of $T(V)$ is then given by $e_se_t=e_{st}$, where $st$ refers to the concatenation of tensor products.

\begin{lem}\label{endocri}
Any map $T(V)\stackrel{\cdot x}{\rightarrow} T(V)$ of right multiplication by a nonzero $x\in T(V)$ is injective. Therefore, every endomorphism of $T(V)$ as a left $T(V)$-module is either zero or injective.
\end{lem}

\pf It suffices to show that the tensor algebra $T(V)$ has no zero-divisors. For basis elements, $e_se_t=e_{st}=e_{s't}=e_{s'}e_t$ implies $st=s't$, thus $s=s'$ by concatenation. Hence $e_s=e_{s'}$. Now since $x\neq0$, we can assume without loss of generality that it has a form $x_0,\sum x_te_t$ or $x_0+\sum x_te_t$ such that the coefficients $x_0$ and $x_t$'s are nonzero. We show for $x=x_0+\sum x_te_t$ and the other cases are similar. Let $y=y_0+\sum y_se_s$ be any element of $T(V)$ for $y_0,y_s\in k$. Then $0=yx=y_0x_0+(\sum y_0x_te_t+\sum y_sx_0e_s)+\sum y_sx_te_{st}$, which implies that $y_0x_0=0$ and $y_sx_t=0$ for all $s,t$. Therefore, $y_0=y_s=0$ for all $s$ and thus $y=0$.~~~\eop

%Let $y\in T(V)$. Suppose $x=u_1+...+u_n$ with each $u_i\neq0$ and $y=v_1+...+v_m$ for $u_i\in V^{\otimes n_i}$ and $v_j\in V^{\otimes m_j}$, where $n_i<n_{i+1}$ and $m_j<m_{j+1}$ for all $i$, $j$. If $yx=0$, then the highest term is $v_mu_n=0$ by the linearly independence of the basis vectors. Therefore, it suffices to show that with respect to the product by the concatenation
%\[
%V^{\otimes s}\times V^{\otimes t}\rightarrow V^{\otimes(s+t)},
%\]
%there are no zero divisors, so that $v_m=0$ and $y=v_1+...+v_{m-1}$. Then the induction implies that $y=0$. Indeed, since the tensor product of vector spaces has a basis given by the tensor product of the basis elements from each vector space, for any two vectors $\a\in V^{\otimes s}$ and $\b\in V^{\otimes t}$, we can assume they have coordinates $\a=a_1e_1+...+a_pe_p$ and $\b=b_1f_1+...+b_qf_q$ with nonzero scalars $a_i,b_j\in k$ and $e_i,f_j$ the basis elements. Then $\a\b=0$ implies that in particular $a_1b_j=0$ for all $j$. Hence $b_j=0$ and $\b=0$, as we required.~~~\eop

%Therefore, $T(V)$ is premonoform by Lemma~\ref{predef} as a left $T(V)$-module.

\begin{lem}\label{containsub}
Suppose the vector space $V$ over a field $k$ has a basis $\{a,b\}$. Let $T(V)/T(V)b$ denote the cokernel of the right multiplication $T(V)\stackrel{\cdot b}{\rightarrow} T(V)$ by $b$. Then the map $T(V)\stackrel{\cdot a}{\rightarrow} T(V)/T(V)b$ of right multiplication by $a$ is injective.
\end{lem}

\pf Let $x=\sum x_se_s$ be an element of $T(V)$ with $x_s\in k$. Suppose $xa+T(V)b=T(V)b$. Then there is $y=\sum y_te_t\in T(V)$ such that
\[
xa-yb=\sum x_s(e_sa)-\sum y_t(e_tb)=0.
\]
Since $e_sa\neq e_tb$ for any $s,t$, it follows that $x_s=y_t=0$ for all $s,t$. Therefore, $x=0$.~~~\eop

%Let $x\in T(V)$, with the coordinate $x=x_1e_1+...+x_ne_n$ using the basis elements $e_i$. Suppose $xa+T(V)b=T(V)b$. Then there is a $y\in T(V)$ with the coordinate $y=y_1f_1+...+y_mf_m$ using the basis elements $f_j$ such that $xa-yb=0$. Then we have
%\[
%x_1e_1a+...+x_ne_na-y_1f_1b-...-y_mf_mb=0,
%\]
%which implies that $x_1=...=x_n=y_1=...=y_m=0$ since $\{e_1a,...,e_na,f_1b,...,f_mb\}$ is part of the basis of $T(V)$ which has no repetition by the concatenation. Therefore, $x=0$.~~~\eop

\begin{prop}\label{prebutnotmo}
The tensor algebra $T(V)$ is premonoform but not monoform as a left $T(V)$-module.
\end{prop}

\pf By Lemma~\ref{predef} and Lemma~\ref{endocri}, the module $T(V)$ is premonoform. By Lemma~\ref{containsub}, the quotient $T(V)/T(V)b$ contains $T(V)$ as a submodule. Thus it is not monoform.~~~\eop

\begin{itemize}
\item \textit{Representation of quivers}
\end{itemize}

In the representation theory of quivers, the indecomposable representations provide many examples of (pre)monoform objects.

We say that the Auslander-Reiten quiver is \textit{directed} if it can be drawn as the AR translations start from the left to the right so that there are no backward arrows. The representations of type $A,D,E$ are of this kind in particular, see Chapter IV.4 in~\cite{quiver}. It follows immediately that any nontrivial (sub)quotients of an indecomposable object, either projective or injective, appear only on its right hand side so that there are no arrows backwards. Therefore, we obtain the following result.

%Let $Q$ be a finite, connected and acyclic quiver, and $A=kQ$ the path algebra of $Q$ over an algebraically close field $k$. Assume $A$ is representation-finite. Then in the category of finite dimensional representations of $A$, an object is monoform if and only if it is premonoform, if and only if it is indecomposable.

\begin{prop}\label{quiver}
Consider the category of finite dimensional representations of a quiver algebra such that the Auslander-Reiten quiver of each representation is connected and directed. Then an object is monoform if and only if it is premonoform, if and only if it is indecomposable.
\end{prop}

%\pf Since the Auslander-Reiten quivers of these representations are connected and directed, .~~~\eop

\begin{itemize}
\item \textit{Modules over a ring}
\end{itemize}

In the category of modules over a ring, a nonzero module is \textit{uniform} if the intersection of any two nonzero submodules is nonzero. Or equivalently, it is an \textit{essential extension} for every nonzero submodule, see Section 1.3F in~\cite{lam}. A uniform module is \textit{strongly uniform} if we replace essential extension by \textit{rational extension}, see Section 1 in~\cite{Storrer}. Note that all of these definitions can be made in an abelian category.

%Comparing with the notions of simple and indecomposable modules, we have the following relationship in the category of modules over a (non)commutative ring with identity.

\begin{prop}~\label{relations}
Let $M$ be an $R$-module over a (non)commutative ring $R$ with an identity. Then for the following concepts, $M$ is (1) simple; (2) strongly uniform; (3) monoform; (4) premonoform; (5) uniform; (6) indecomposable, we have implications
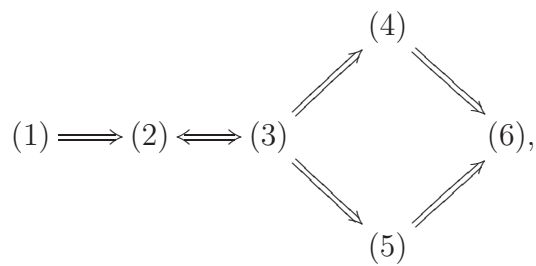
\[
\xymatrix{
&&&(4) \ar@{=>}[dr] &\\
(1)\ar@{=>}[r] &(2) \ar@{<=>}[r] &(3)\ar@{=>}[dr] \ar@{=>}[ur] &&(6),\\
&&& (5)\ar@{=>}[ur]
}
\]
where (4) and (5) are not comparable.
\end{prop}

%$(1)\Rightarrow(2)\Leftrightarrow (3)\Rightarrow \{(4), (5)\}\Rightarrow (6)$,

%In fact, any nonrational extension of a module $M$ gives a nonzero map $f:A/N\rightarrow M$ for some $N<A<M$, whose coimage is a subobject of $M$ and any nontrivial subquotient $A/N$ of $M$ gives a nonzero map

\pf The equivalence $(2)\Leftrightarrow(3)$ holds by Proposition 3.8.6 in~\cite{lam}. The implication $(1)\Rightarrow(3)$ holds because the simplicity of $M$ gives no nontrivial subquotient of $M$. Also, $(2)\Rightarrow (5)$ and $(3)\Rightarrow (4)$ are true by the definition. For $(4)\Rightarrow (6)$, any decomposition $M=M_1\oplus M_2$ gives an injection $M/M_1\cong M_2\hookrightarrow M$. Since any decomposition $M=M_1\oplus M_2$ gives a trivial intersection $M_1\cap M_2=0$, the implication $(5)\Rightarrow (6)$ holds. The incomparable relation between (4) and (5) is demonstrated in Example~\ref{compareex}.~~~\eop

The following examples show that the implications in Proposition~\ref{relations} are all strict, and that (4), (5) are not comparable. Some of them can be found in~\cite{lam}.

\begin{ex}\label{compareex}
(1) Any simple $R$-module is of the form $R/\mm$ for some maximal ideal $\mm$. The quotient $R/(x)$ of the polynomial ring $R=\QQ[x,y]$ gives a monoform $R$-module that is not simple. Indeed, $R/(x)$ has the quotient ideal $(x,y)/(x)$ as a nontrivial submodule.

\noindent(2) The cyclic group $\ZZ/4$ is uniform but neither premonoform nor monoform. More generally, we can consider the cyclic groups $\ZZ/p^n$ of order $p^n$ for $n>1$.

\noindent(3) The tensor algebra $T(V)$ of a $k$-vector space $V$ with basis $\{a,b\}$ is a premonoform $T(V)$-module by Proposition~\ref{prebutnotmo} but not uniform. Indeed, we have $T(V)a\cap T(V)b=0$, thanks to the form of the basis elements given by the monomials in $a,b$. For example, $T(V)$ has $\{a^2,ba,ab,b^2\}$ as the monomial basis elements of degree 2. Also, notice that this gives another proof of Proposition~\ref{prebutnotmo}.

\noindent(4) The commutative $\QQ$-algebra $R=\QQ[x,y]/(x^2,xy,y^2)$ is indecomposable as an $R$-module but not uniform since in $R$ there is a trivial intersection $R x\cap R y=0$. The algebra $R$ is not premonoform either since there is a nonzero $R$-module homomorphism $f:R\rightarrow R$ defined by $f(1)=x$. However, $f(y)=0$, so $f$ is not injective.

%\noindent(4) The commutative $\QQ$-algebra $R=\QQ[x_1,x_2,...]$ with countably variables substituting to the relations $x_1^2=0$ and $x_n^2=x_{n-1}$ for $n>1$ is uniform but not strongly uniform
\end{ex}

\begin{itemize}
\item \textit{Finitely generated modules over a commutative noetherian ring with identity}
\end{itemize}

\begin{prop}\label{precoincide}
In the category of finitely generated modules over a commutative noetherian ring $R$ with an identity, an $R$-module $M$ is premonoform if and only if it is monoform.
\end{prop}

\pf As always monoform implies premonoform. Now let $M$ be premonoform and suppose it is not monoform. Then by the definition there exists an $H\neq0$ as a common submodule of both $M$ and $X=M/N$ for some nonzero submodule $N\subseteq M$. Thus there is an associated prime ideal $\mp\in\Ass H$ such that $R/\mp\hookrightarrow H$. In particular, $X_\mp\neq0$ since $\mp\in\Ass H\subseteq\Ass X\subseteq\Supp X$. Passing to the local case at $\mp$, we obtain on the one hand that $X_\mp\otimes R/\mp\cong X_\mp\otimes R_\mp\otimes R/\mp\cong X_\mp\otimes k(\mp)\cong \bigoplus k(\mp)$ since $X_\mp$ is finitely generated, and on the other hand $X_\mp\otimes R/\mp\cong X_\mp/\mp X_\mp$ which is nonzero by the Nakayama Lemma. Therefore, we have an exact sequence
\[
0\rightarrow \mp X_\mp\rightarrow X_\mp\stackrel{f}{\rightarrow} \bigoplus k(\mp)\rightarrow 0,
\]
in which $f$ is nonzero. Moreover, composed with a projection onto one copy $k(\mp)\cong (R/\mp)_\mp$, by Proposition 2.10 in~\cite{eisen} the map $f$ can be lifted into a nonzero map $g:X\rightarrow R/\mp$ since $X$ is finitely generated. Hence there is a nonzero map $M/N=X\stackrel{g}{\rightarrow} R/\mp\hookrightarrow H\subseteq M$, a contradiction to $M$ being premonoform by Lemma~\ref{predef}.~~~\eop

\section{Spectrum, support and topology}\label{spectrum}

%From the way it defines monoform objects and as it shows in~\cite{Kanda} that the support can be used to classify Serre subcategories, similarly, premonoform objects classify nullity classes.

%those full subcategories closed under quotients and extensions. They are also identified with torsion classes in some context.

In this section, we will define our spectrum together with the support of objects and subcategories. These are the basic concepts we need to establish the classification theorem.

Let $\al$ be an abelian category. Denote by $\Spec_0\al$ the collection of isomorphism classes of premonoform objects in $\al$. Let $A\in\al$. Denote by $\overline{C(A)}$ the smallest nullity class containing $A$, that is, the intersection of all nullity classes containing $A$. We also call $\overline{C(A)}$ the nullity class generated by $A$.

\begin{de}
Define $A\sim B$ if and only if $\overline{C(A)}=\overline{C(B)}$, i.e. they generate the same nullity class.
\end{de}

%We give another equivalent description later by using support given as Proposition~\ref{equiv}.\\

%\Supp M=\{\overline{H}\in\Spec\al~|~\exists H'\in\overline{H}~s.t. ~H'~\text{is a quotient of }M\}.
%%%%%%%%%%%%%%%%%%%%%%%%%%%%%%%%%%%%%%%%%%%%%%
%this also gives a topological space on $\Spec\al$ but requires that the category $\al$ is Noetherian.

It is clear that $\sim$ is an equivalence relation on $\Spec_0\al$, and we denote by $\Spec\al=\Spec_0\al/\sim$ the collection of the equivalence classes $[H]$ of premonoform objects, called the \textit{spectrum} of $\al$ which is either a class or a set. A sufficient condition for the spectrum $\Spec\al$ being a set is that the category $\al$ is \textit{essentially small}, i.e. the collection of the isomorphism classes of objects in $\al$ forms a set. For any $M\in \al$, define the \textit{support of an object} as
\[
\Supp M=\{[H]\in\Spec\al~|~\exists~ H'\in[H]~s.t. ~H'\in\overline{C(M)}\}=\{[H]\in\Spec\al~|~
[H]\subseteq\overline{C(M)}\}.
\]
The \textit{support of a subcategory} $\nl\subseteq\al$ is defined as the union
\[
\Supp~\nl=\bigcup_{M\in \nl}\Supp~M.
\]
We can also define
\[
\Supp^{-1}\Phi=\{M\in\al~|~\Supp M\subseteq\Phi\}
\]
as a full subcategory of $\al$ for any subclass $\Phi\subseteq\Spec\al$. However, this is not always a nullity class for an arbitrary subclass $\Phi$ but only for those closed and extension closed ones, as we will see in Theorem~\ref{gensttment}.

These definitions are similar to Kanda's in~\cite{Kanda}, in which the \textit{atom spectrum} $\ASpec\al$ of an abelian category $\al$ consists of equivalence classes of monoform objects, and the \textit{atom support} of an object $M$ is given by $\ASupp M=\{[H]\in\ASpec\al~|~\exists~ H'\in[H]~s.t. ~H'~\text{is a subquotient of}~M\}$. Here for two monoform objects, $H\sim H'$ are \textit{atom equivalent} if they share a common nonzero subobject.

\begin{rem}
Consider the category of abelian groups and the objects $\QQ$ and $\ZZ$. They are clearly monoform and thus premonoform, because no torsion-free groups contain a torsion subgroup. With respect to the atom equivalent, we have $\QQ\sim\ZZ$. However, they are not equivalent in our sense since $\QQ\notin\overline{C(\ZZ)}=\langle\langle\ZZ\rangle_{\quot}\rangle_{\ext}$.
\end{rem}

%Worth mentioning that such defined concept of support shares the feature

%there is $H\in\overline{M}$ such that $\Supp H\subseteq\Phi$.

\begin{de}
Let $\al$ be an abelian category and $\Spec\al$ its spectrum (either a set or a class), $\Phi\subseteq\Spec\al$ a subclass. The subclass $\Phi$ is \textit{closed} if for every $[M]\in \Phi$, $\Supp M\subseteq \Phi$. The subclass $\Phi$ is \textit{extension closed} if whenever there is a short exact sequence
\[
0\rightarrow M\rightarrow X\rightarrow N\rightarrow 0,
\]
if $\Supp M,\Supp N\subseteq\Phi$ then $\Supp X\subseteq \Phi$.
\end{de}

%There are some immediate consequences as the next proposition and lemma.

%\begin{prop}
%If $\al$ is an abelian category such that $\Spec\al$ forms a set. Then the open sets satisfy the axiom of topology.
%\end{prop}

%\pf Suppose $\Phi,\Psi\subseteq \Spec\al$ are open and $\overline{M}\in\Phi\cap\Psi$. Then $\Supp M\subseteq\Phi\cap\Psi$ hence $\Phi\cap\Psi$ is open. Also, if $M\in\bigcup\Phi_i$ for $\Phi_i$ open, then $M\in\Phi_i$ for some $i$ implies $\Supp M\subseteq\Phi_i\subseteq\bigcup \Phi_i$ thus $\bigcup \Phi_i$ remains open.~~~\eop

\begin{lem}\label{small properties}
Let $M,N$ be objects in an abelian category $\al$. Then

\noindent (1) if $M$ is premonoform, then $[M]\in\Supp M$;

\noindent(2) if $M$ is premonoform, then $\overline{C(M)}\subseteq \overline{C(N)}$ if and only if $\Supp M\subseteq \Supp N$. In particular, if both $M$ and $N$ are premonoform, then $M\sim N$ if and only if $\Supp M=\Supp N$;

\noindent(3) if $M\rightarrow N\rightarrow 0$ is exact, then $\Supp N\subseteq\Supp M$;

\noindent(4) the subclass $\Supp M$ is closed.

%\noindent(4) Assume $\Spec\al$ is a set. Then $\{\Supp M \}_{\overline{M}\in\Spec\al}$ forms a topology base for $\Spec\al$.

%\noindent (5) $M$ premonoform implies $|\Ass M|=1$.
\end{lem}

%For every $\overline{A}\in\Supp N$, there is $H\in\overline{A}$ such that $H$ is a quotient of $N$, hence a quotient of $M$, i.e. $\overline{A}\in\Supp M$. This proves (2).

%$\overline{M}\in \Supp M\subseteq\Supp N$ implies that there is an $M'\sim M$ such that $M'\in \overline{C(N)}$. So $\overline{C(M)}=\overline{C(M')}\subseteq \overline{C(N)}$.

\pf (1) Since $M$ is premonoform and $M\in\overline{C(M)}$, thus $[M]\in\Supp M$ by the definition. For (2), $\Supp M\subseteq\Supp N$ implies that in particular $[M]\in\Supp N$ by (1) and thus $\overline{C(M)}\subseteq\overline{C(N)}$. The converse is true by the definition. Property (3) holds because nullity class is closed under quotients. For (4), since $[H]\in\Supp M$ implies $H\in\overline{C(M)}$, thus $[H]\in\Supp H\subseteq\Supp M$ by (1) and (2).~~~\eop

%take any $\overline{M'}\in\Supp M$. Then there is an $H\in\overline{M'}$ such that $H\in \overline{C(M)}$, thus $\Supp M'=\Supp H\subseteq \Supp M$ by (2).~~~\eop

%If $H$ is a quotient of $N$ which is a quotient of $M$ then $H$ is a also a quotient of $M$. This gives (2) and (3). Clearly, this collection of open sets cover $\Spec\al$ by (1). Now if $H\in\Supp M\cap\Supp N$ for $M,N\in\Spec\al$ then $H\in\Supp H\subseteq\Supp M\cap\Supp N$ by (1) and (2).~~~\eop

\begin{lem}
A subclass $\Phi\subseteq\Spec\al$ is closed if and only if for every $[M]\in\Phi$ there is an $H\in\al$ such that $[M]\in\Supp H\subseteq \Phi$.
\end{lem}

%Clearly, the openness of $\Phi$ implies even more that such $H$ can be chosen in the equivalence class $\overline{M}$.

\pf The necessity holds by the definition. For the sufficiency, let $[M]\in\Phi$ and suppose $[M]\in\Supp H\subseteq \Phi$ for some $H\in\al$. Then $[M]\subseteq\overline{C(H)}$, hence $\Supp M\subseteq \Supp H\subseteq\Phi$ by Lemma~\ref{small properties}, as required.~~~\eop

\begin{prop}\label{topspace}
Assume $\al$ is an abelian category such that $\Spec\al$ forms a set. Then the collection of the closed subsets of $\Spec\al$ indeed forms a topology of closed subsets. This topology is specialization closed.
\end{prop}

\pf It is clear by the definition that the empty set and the whole set $\Spec\al$ are closed. Let $\Phi_i$ be a closed subset for every $i\in I$. Suppose $[M]\in\cup_{i\in I}\Phi_i$. Then $[M]\in \Phi_i$ for some $i$ implies that $\Supp M\subseteq\Phi_i\subseteq\cup_{i\in I}\Phi_i$. Hence $\cup_{i\in I}\Phi_i$ is closed. Now suppose $[M]\in\cap_{i\in I}\Phi_i$. Then $[M]\in \Phi_i$ implies that $\Supp M\subseteq\Phi_i$, for every $i\in I$. Hence $\Supp M\subseteq\cap_{i\in I}\Phi_i$, that is, $\cap_{i\in I}\Phi_i$ is also closed.~~~\eop

We denote by $\Spec_{\noeth}\al$ the subset (or subclass) of $\Spec\al$ in which every equivalence class $[A]\in\Spec_{\noeth}\al$ has a noetherian representative $A\in\al$, then the subset $\Spec_\noeth\al$ is shown to be a topological subspace of $\Spec\al$ by the same argument of Proposition~\ref{topspace}, noticing that every closed subset of $\Spec_\noeth\al$ is also closed in $\Spec\al$.

\section{Classification of nullity classes}\label{maintheorem}

Kanda shows in~\cite{Kanda} that in a noetherian abelian category $\al$ there is a bijective correspondence from the collection of open subsets of the atom spectrum $\ASpec\al$ of equivalence classes of monoform objects, to the collection of Serre subcategories of $\al$. In this section, enlightened by this idea we give a classification of nullity classes via our spectrum $\Spec\al$ of equivalence classes of premonform objects in $\al$.

%We characterize the premonoform object in a slightly different way due to its construction, in order to give a classification of the nullity class. This classification resembles that of the Serre subcategories in a noetherian abelian category shown by Kanda in~\cite{Kanda}.

The following key lemma gives a necessary and  sufficient condition for an object to be premonoform.

\begin{lem}\label{notprem}
Let $\al$ be an abelian category and $M\in\al$. To say that $M$ is not premonomoform is equivalent to saying that $M$ lies in the nullity class $\nl$ generated by all quotients $M/N$ with nonzero subobjects $N\leq M$.
\end{lem}

%this is also true for torsion classes, noticing that $M$ is not in the torsion class iff there are no maps from $M/N$ to $M$ for every $N\neq0$.

\pf Observe that for the case $M=0$, the statement holds obviously. So assume from now on $M\neq0$. If $M$ is not premonoform, then $M$ has a nontrivial quotient $M/N$ as its subobject. That is, $M/N\cong N'\leq M$ for some $N'$. Then the short exact sequence
\[
0\rightarrow M/N\rightarrow M\rightarrow M/N'\rightarrow 0
\]
implies that $M\in\nl$ since both $M/N,M/N'\in\nl$ and $\nl$ is closed under extensions. Conversely, suppose $M\in\nl$. Then by Lemma~\ref{ext}, it is an iterated extension of quotients of $M$. So suppose $M$ is an $n$-extension but not an $(n-1)$-extension for $n\geq1$, see Notation~\ref{next}. Thus there is a short exact sequence
\[
0\rightarrow M/N\rightarrow M\rightarrow M'\rightarrow 0
\]
such that both $M/N$ and $M'$ are nonzero, where $M'$ is an $(n-1)$-extension of the quotients of $M$. Hence $M$ has a nontrivial quotient $M/N$ as a subobject. Therefore, $M$ is not premonoform.~~~\eop

%Our result still holds as long as every object is Noetherian.

%We say that an object in an abelian category is \textit{Noetherian} if whenever there is an increasing chain of subobjects it becomes stationary after finite many steps.

Now we are ready to give a criterion using support to detect when an object is in a nullity class.

\begin{prop}\label{criterion}
Let $\al$ be an abelian category and $\nl$ a nullity class, $M\in\al$ a noetherian object. Then $M\in \nl$ if and only if $\Supp M\subseteq \Supp \nl$.

%(2) If $\al$ is well-powered and subcomplete, then $M\in\nl$ if and only if $\Supp M\subseteq \Supp\nl$.
\end{prop}

\pf One implication is obtained by the definition of support without any extra assumption. Indeed, $M\in\nl$ implies $\Supp M\subseteq\bigcup_{N\in\nl}\Supp N=\Supp\nl$. Now suppose $\Supp M\subseteq \Supp \nl$ but $M\notin\nl$. Then consider the collection $\sl=\{N\leq M~|~M/N\notin\nl\}$ of subobjects of $M$. Notice that $\sl$ is nonempty since $M\notin\nl$ implies that $N=0$ is one such element. Hence $\sl$ contains a maximal subobject, say $N_0\leq M$ such that $M/N_0\notin\nl$ since $M$ is noetherian. Then we claim that $M/N_0$ is premonoform by Lemma~\ref{notprem} since otherwise $M/N_0$ must lie in the nullity class generated by the nontrivial quotients of $M/N_0$ all belonging to $\nl$ by the maximality of $N_0$, and thus $M/N_0\in\nl$ as well, a contradiction. So $[M/N_0]\in\Supp M\subseteq\Supp\nl$. Hence there is an $X\in\nl$ such that $[M/N_0]\in\Supp X$, that is, $[M/N_0]\subseteq\overline{C(X)}$. Therefore,
\[
M/N_0\in[M/N_0]\subseteq\overline{C(X)}\subseteq\nl,
\]
a contradiction.~~~\eop

We will see from Example~\ref{suppnotdet} that noetherian is a required hypothesis in Proposition~\ref{criterion}.

\begin{nota}\label{notationforc}
Let $\al$ be an abelian category.

\noindent(1) Denote by $\nl$ a nullity class and $\tl$ a torsion class in $\al$. For example, Proposition~\ref{torsionnullity} says that if $\al$ is noetherian, then the two collections $\{\nl\subseteq\al\}=\{\tl\subseteq\al\}$ coincide.

\noindent(2) Denote by $\nl_{\noeth}$ a nullity class consisting of noetherian objects.

%and $\nl_{\prem}$ a nullity class generated by the premonoform objects, i.e. it is the intersection of all nullity classes containing every premonoform object $M$ in $\nl$. Similarly, we can define for torsion classes.

\noindent(3) Denote by $S^{\ext}_c$ a subclass or a subset of $\Spec\al$ (or $\Spec_{\noeth}\al$) that is closed and extension closed.

%\noindent(4) Denote by $\Spec_{\noeth}\al$ the subclass of $\Spec\al$ in which every equivalence class $[A]\in\Spec_{\noeth}\al$ has a noetherian representative $A\in\al$.
\end{nota}

The following statement is an easy generalization of Proposition 6.3 in~\cite{Atiyah}, which is stated for the noetherian modules.

\begin{lem}\label{noethobj}
Let $\al$ be an abelian category. Suppose there is a short exact sequence
\[
0\rightarrow N'\stackrel{f}{\rightarrow} M\stackrel{g}{\rightarrow} N\rightarrow 0
\]
in $\al$. Then $M$ is noetherian if and only if both $N, N'$ are noetherian. This implies for two equivalent premonoform objects $A\sim B$, $A$ is noetherian if and only if $B$ is noetherian.
\end{lem}

\pf Since any chain of subjects in $N'$ or $N$ also gives a chain of subobjects in $M$, the necessity follows immediately. Conversely, suppose $\{M_i\}$ is a chain in $M$. Then $\{f^{-1}M_i\}$ and $\{g(M_i)\}$ are chains in $N',N$ respectively, so that for large enough index $i$ that the stationary of both chains $\{f^{-1}M_i\},\{g(M_i)\}$ implies that of the chain $\{M_i\}$. The last statement then follows immediate from Lemma~\ref{ext}.~~~\eop

\begin{prop}\label{gennull}
There is an invariant of nullity classes in an abelian category $\al$ given by the support
\[
\Supp:\{\nl\subseteq\al\}\rightarrow\{S_c\subseteq\Spec\al\},
\]
where $S_c$ represents a closed subclass or subset of $\Spec\al$.
\end{prop}

\pf For any nullity class $\nl$, the support $\Supp\nl$ is closed since for every $[H]\in\Supp \nl$ there is an $X\in\nl$ such that $[H]\in\Supp X\subseteq\Supp\nl$. Also, for any two nullity classes $\nl_1,\nl_2$, if $\Supp \nl_1\neq\Supp\nl_2$, say there is $[M]\in\Supp\nl_1-\Supp\nl_2$, then $[M]\subseteq\nl_1$ and $[M]\nsubseteq\nl_2$. Hence $M\in\nl_1-\nl_2$ and thus $\nl_1\neq\nl_2$.~~~\eop

\begin{ex}\label{zpinfty}
In the category of abelian groups, consider the the nullity class $\nl$ generated by $\ZZ_{(p)}$ and its support $\Phi=\Supp\nl$. First observe that $\Supp\ZZ/p^\infty=\emptyset$. Indeed, $\ZZ/p^\infty$ is not premonoform due to the short exact sequence
\[
0\rightarrow\ZZ/p\rightarrow\ZZ/p^\infty\stackrel{\cdot p}{\rightarrow}\ZZ/p^\infty\rightarrow0.
\]
Since $\overline{C(\ZZ/p^\infty)}=\langle\langle\ZZ/p^\infty\rangle_{\quot}\rangle_{\ext}=\langle\ZZ/p^\infty\rangle_{\ext}$ by Lemma~\ref{ext} and $\ZZ/p^\infty$ is injective, every object in $\overline{C(\ZZ/p^\infty)}$ is a direct sum of finite copies of $\ZZ/p^\infty$. Hence it has empty support. Also notice that there is a short exact sequence
\[
0\rightarrow \ZZ_{(p)}\hookrightarrow\QQ\rightarrow\ZZ/p^\infty\rightarrow0,
\]
such that both $\Supp \ZZ_{(p)}$ and $\Supp\ZZ/p^\infty$ are contained in $\Phi$. However, $[\QQ]\in\Supp\QQ-\Phi$. This implies that the support $\Phi$ of $\nl$ is not extension closed.
\end{ex}

The next lemma holds for an arbitrary abelian category.

\begin{lem}\label{lemsupp}
Let $\al$ be an abelian category and $\Phi$ a closed and extension closed subclass or subset of $\Spec\al$. Then

\noindent(1) $\Supp^{-1}\Phi$ is a nullity class;

\noindent(2) $\Supp\Supp^{-1}\Phi=\Phi$.
\end{lem}

\pf For (1), the category $\Supp^{-1}\Phi$ is closed under extensions since whenever $M,N\in\Supp^{-1}\Phi$, the support satisfies $\Supp X\subseteq\Phi$ for any extension $X$ of $N$ by $M$, because $\Phi$ is extension closed. Thus $X\in\Supp^{-1}\Phi$. Also, $\Supp^{-1}\Phi$ is closed under quotients by Lemma~\ref{small properties}. Therefore, $\Supp^{-1}\Phi$ is a nullity class. For (2), suppose $[M]\in\Supp\Supp^{-1}\Phi$. That is, $[M]\subseteq \overline{C(X)}$ for some $X\in\Supp^{-1}\Phi$. It follows that $\Supp M\subseteq\Supp X\subseteq \Phi$ by Lemma~\ref{small properties} and thus $M\in\Supp^{-1}\Phi$. Therefore, $[M]\in\Supp M\subseteq\Phi$ and thus $\Supp\Supp^{-1}\Phi\subseteq\Phi$. Conversely, suppose $[M]\in\Phi$. Since $\Phi$ is closed, there is an $H\in\al$ such that $[M]\in\Supp H\subseteq\Phi$. In other words, $[M]\in\Supp\Supp^{-1}\Phi$. Therefore, $\Phi\subseteq\Supp\Supp^{-1}\Phi$, as required.~~~\eop

We can now prove our main result.

\begin{thm}\label{gensttment}
Let $\al$ be an abelian category. Then the nullity classes of noetherian objects are classified by the closed and extension closed subclasses of the spectrum $\Spec_{\noeth}\al$. More precisely, there is an order preserving bijection
\[
\Supp:\{\nl_{\noeth}\subseteq\al\}\stackrel{\sim}{\rightleftarrows}\{S^{\ext}_c\subseteq\Spec_{\noeth}\al\}: \Supp^{-1}
\]
in which $\Supp^{-1}\Phi=\{M\in\al~|M~\text{is noetherian and}~\Supp M\in\Phi\}$.
\end{thm}

\pf We will first see that the support of a nullity class is extension closed. So assume $\nl$ is a nullity class of noetherian objects in $\al$ and $M,N\in \al$ are noetherian such that $\Supp M,\Supp N\subseteq\Supp \nl$. In particular, $M,N\in \nl$ by Proposition~\ref{criterion}. Suppose $X$ is an extension
\[
0\rightarrow M\rightarrow X\rightarrow N\rightarrow 0
\]
of $N$ by $M$. Then $X$ is noetherian by Lemma~\ref{noethobj} and thus $X\in\nl$ since $\nl$ is closed under extensions. Hence $\Supp X\subseteq\Supp \nl$ and $\Supp\nl$ is extension closed. Furthermore, $\Supp\nl$ is closed by Proposition~\ref{gennull}. On the other hand, given any closed and extension closed subclass $\Phi$ of $\Spec_{\noeth}\al$, $\Supp^{-1}\Phi$ is a nullity class by Lemma~\ref{lemsupp}. Hence $\Supp^{-1}\Phi$ is a nullity class of noetherian objects by Lemma~\ref{noethobj}, again. Now suppose $\nl$ is a nullity class of noetherian objects. Then since $\Supp M\subseteq\Supp \nl$ if and only if $M\in\nl$ by Proposition~\ref{criterion}, $\Supp^{-1}\Supp\nl=\nl$. Furthermore, for a closed and extension closed $\Phi$, we have $\Supp\Supp^{-1}\Phi=\Phi$ by Lemma~\ref{lemsupp}. This completes the proof.~~~\eop

We will see from Example~\ref{a2} in the following section that the extension closed condition is necessary. We also wonder what categories do not require this condition.

%Now suppose $\nl$ is a nullity class. Then since $\Supp M\subseteq\Supp \nl$ if and only if $M\in\nl$ by Proposition~\ref{criterion}, hence $\Supp^{-1}\Supp\nl=\nl$. On the other hand, suppose $\overline{H}\in\Supp\Supp^{-1}\Phi$, i.e. there is $H'\in\overline{H}$ such that it is a quotient of $X\in\Supp^{-1}\Phi$. Hence $\Supp H'\subseteq\Supp X\subseteq \Phi$ so that $H'\in\Supp^{-1}\Phi$. Therefore, $H\in\overline{C(H)}=\overline{C(H')}\subseteq\Supp^{-1}\Phi$, i.e. $\overline{H}\in\Supp H\subseteq\Phi$. Now suppose $\overline{H}\in\Phi$. Then there is $M\in\al$ such that $\overline{H}\in\Supp M\subseteq\Phi$ by openness. That is, there is $H'\in\overline{H}$ such that it is a quotient of $M\in\Supp^{-1}\Phi$. Therefore, $\overline{H}\in\Supp \Supp^{-1}\Phi$.~~~\eop

\section{Classification of torsion classes}\label{variants}

In this section, we give an immediate consequence of Theorem~\ref{gensttment}. Notice that the spectrum $\Spec\al$ becomes a set if the abelian category $\al$ is noetherian, and thus $\nl=\nl_{\mathrm{noeth}}$, $\Spec\al=\Spec_{\noeth}\al$ and nullity classes are torsion classes by Proposition~\ref{torsionnullity}. See Notation~\ref{notationforc}.

\begin{thm}\label{bij}
Suppose $\al$ is a noetherian abelian category. Then there is an order preserving bijection
\[
\Supp:\{\tl\subseteq\al\}\stackrel{\sim}{\rightleftarrows}\{S^{\ext}_c\subseteq\Spec\al\}:\Supp^{-1}.
\]
In other words, the torsion classes of $\al$ are classified by the closed and extension closed subsets of the spectrum $\Spec \al$.
\end{thm}

\begin{rem}
Since the collection of nullity classes in the category of abelian groups does not form a set, by Proposition 4.7 in~\cite{Stanw2} (see Corollary 8.4 in~\cite{Stan} for a derived category version), we cannot obtain a classification via a topological space.

\end{rem}

We end this paper with two examples.

\begin{ex}\label{suppnotdet}
Consider the category of abelian groups and the nullity class $\overline{C(\ZZ/p^\infty)}$ generated by the non-noetherian abelian group  $\ZZ/p^\infty$. Since it has empty support by Example~\ref{zpinfty}, our system of support cannot detect this nullity class $\overline{C(\ZZ/p^\infty)}$. The same problem happens for torsion classes, by using a similar notion of spectrum and support. Denote by $T(\ZZ/p^\infty)$ the torsion class generated by $\ZZ/p^\infty$, which is a cocomplete nullity class by Proposition~\ref{torsionnullity}. Next suppose $M$ is a nontrivial abelian group and there is a surjection $f:\oplus\ZZ/p^\infty\twoheadrightarrow M$. Since $\ZZ$ is noetherian, the direct sum $\oplus\ZZ/p^\infty$ is injective by Theorem 3.46 in Lam~\cite{lam}, thus it is divisible by Theorem 7.1 in Hilton and Stammbach~\cite{hilt}  and so is $M$ by Proposition 7.2 also in~\cite{hilt}. Therefore, $M$ is injective and it is a direct sum of copies of $\ZZ/q^\infty$ for some primes $q$. Since $\Hom_\ZZ(\oplus\ZZ/p^\infty,\oplus_q\ZZ/q^\infty)
\cong\prod\Hom_\ZZ(\ZZ/p^\infty,\oplus_q\ZZ/q^\infty)$, the map $f$ is determined by its components $f_i:\ZZ/p^\infty\rightarrow\oplus_q\ZZ/q^\infty$, say $f_i(1)=(x_j)_j$ for only finitely many nonzero $x_j\in\ZZ/q^\infty$. However, by an argument of element orders we see that $f_i=0$ unless $q=p$ for all $q$. Hence $M\cong\oplus\ZZ/p^\infty$, as needed.

%Thus thanks to the previous argument, it suffices to show that for $M\in T(\ZZ/p^\infty)$ any surjection $f:\oplus\ZZ/p^\infty\rightarrow M$ is an isomorphism, so that every object in $T(\ZZ/p^\infty)$ is an arbitrary direct sum of copies of $\ZZ/p^\infty$. Indeed, each direct summand of $\oplus\ZZ/p^\infty$ is a summand of $M$
%\[
%\xymatrix{
%&\ZZ/p^\infty\ar@{->>}[d]\ar@{^(->}[r]&\oplus\ZZ/p^\infty\ar^f@{->>}[r] &M\\
%\ZZ/p^\infty\ar^{\cong}[r]&\ZZ/p^\infty/\Ker(f)\ar@{^(->}[urr]&&
%}
%\]
%by the above short exact sequence and the injectivity of $\ZZ/p^\infty$. Therefore, to get rid of the noetherian condition, there are some other things to consider.
\end{ex}

%Need to show that for $M$ premonoform, $M\in C(\Phi)$ iff $\overline{M}\in\Phi$.

%We end this section with one example.

%We compute the spectrum of the category $kA_2\text{-}\rep$ consisting of finite dimensional representations of the quiver $A_2:\stackrel{1}{\circ}\leftarrow\stackrel{2}{\circ}$ over a field as follows.

\begin{ex}\label{a2}
Let $\al$ be the category of finite dimensional representations of the quiver $A_2:\stackrel{1}{\circ}\leftarrow\stackrel{2}{\circ}$ over a field $k$. Let $a=P_1,b=P_2,c=S_2$ denote the three indecomposables.

\noindent(1) The indecomposables give all the premonoform objects by Proposition~\ref{quiver}. Since no pair of them are equivalent, the underlying set of the spectrum is
\[
\Spec \al=\{a,b,c\}.
\]
whose topology of closed subsets other than the empty set and $\Spec\al$ is given by $\Supp~a=\{a\}$, $\Supp~c=\{c\}$, $\Supp~b=\{b,c\}$ and the union $\{a,c\}$. Except $\{a,c\}$, each closed subset is extension closed, thus corresponds to a nullity class uniquely by Theorem~\ref{bij}, for instance.

\noindent(2) The computation in (1) also implies that our support, unlike the atom support, does not respect extensions since $\Supp ~b\nsubseteq\Supp~ a\cup \Supp ~c$.

\noindent(3) As illustrated in the diagram of the lattice $\Phi$ of nullity classes in $\al$
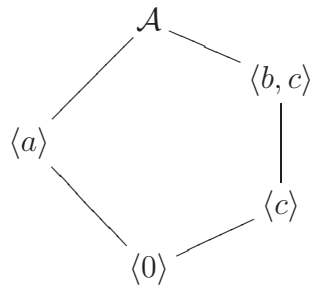
\[
\xymatrix@R=0.5em{
&\al\ar@{-}[dr]&\\
&&\langle b,c\rangle\ar@{-}[dd]\\
\langle a\rangle\ar@{-}[uur]\ar@{-}[ddr]&&\\
&&\langle c\rangle\\
&\langle0\rangle\ar@{-}[ur]&
}
\]
the lattice $\Phi$ is not distributive by Theorem 4.10 (ii) in~\cite{DaveyPriestley} since it is isomorphic to a pentagon, so that it cannot be isomorphic to a lattice of closed subsets of a topological space.

Therefore, the condition of being extension-closed in our main result is necessary, if we use a topological space to classify the nullity classes by establishing a bijective correspondence from the collection of closed subsets to that of nullity classes.
\end{ex}

%\begin{ex}
%By reviewing the example of the tensor algebra $T(V)$ with the $k$-vector space $V$ of basis $\{ a,b\}$ in Section~\ref{relatedcon}, we consider the nullity class $\overline{C(T(V))}$ generated by the premonoform object $T(V)$. It is not a nullity class of noetherian objects since $T(V)$ is not noetherian. Therefore, Theorem~\ref{nullgenbyprem} is indeed a generalization of Theorem~\ref{bij}.

%Also, there are many nullity classes not generated by the premonoforms either, such as the nullity classes $\overline{C(W)}$ of vector spaces generated by an infinite dimensional vector space $W$ with a basis of cardinality say $\aleph_0$, in which the objects are those vector spaces with a basis of cardinality $\leq\aleph_0$.
%\end{ex}

%There are still many interesting questions to ask, for example, what can we say about torsion classes and localizing subcategories?

\bibliography{mybib}
\bibliographystyle{abbrv}

\end{document}